\documentclass[11pt,a4paper]{article}
\usepackage{latexsym}
\usepackage{amsmath}
\usepackage{amssymb}
\usepackage{amsthm}

\usepackage[pagebackref, plainpages=false, pdfpagelabels]{hyperref}

\swapnumbers

\theoremstyle{plain}

\newtheorem{thm}[subsection]{Theorem}
\newtheorem{prop}[subsection]{Proposition}
\newtheorem{lemma}[subsection]{Lemma}
\newtheorem{cor}[subsection]{Corollary}

\theoremstyle{definition}
\newtheorem{remark}[subsection]{Remark}

\newcommand{\abs}[1]{\lvert#1\rvert}        

\newcommand{\form}[2]{\langle #1 | #2 \rangle}	

\newcommand{\LemRef}[1]{Lemma~\ref{#1}}
\newcommand{\PropRef}[1]{Proposition~\ref{#1}}
\newcommand{\ThmRef}[1]{Theorem~\ref{#1}}

\newcommand {\Aff} {\mathrm{Aff}}

\newcommand {\Fi} [1] {{\rm Fi}_{#1}}

\newcommand {\A} {\mathcal{A}}
\newcommand {\V} {\mathcal{V}}
\newcommand {\I} {\mathcal{I}}
\newcommand {\J} {\mathcal{J}}
\newcommand {\IAff} {\I_\mathrm{Aff}}
\newcommand {\OA} {\bar{\mathcal{A}}}

\newcommand {\F} {\mathbb{F}}

\newcommand{\SU}{{\rm SU}}

\renewcommand{\O}{{\rm O}}
\newcommand{\PO}{{\rm PO}}
\newcommand{\PSU}{{\rm PSU}}
\newcommand{\POmega}{{\rm P}\Omega}
\newcommand{\Sym}{{\rm Sym}}
\renewcommand{\split}{\colon\!}

\usepackage{color}      

\begin{document}
\title {Lie algebras and $3$-transpositions}
\author{H. Cuypers, M. Horn, J. in 't panhuis, S. Shpectorov}
\date{\today}
\maketitle


\begin{abstract}
\noindent
We describe a construction of an algebra over the field of order $2$ starting from a conjugacy class
 of  $3$-transpositions in a group. In particular, we determine which  simple Lie algebras arise by this 
construction. Among other things, this construction yields a natural embedding of the sporadic simple 
group $\Fi{22}$ in the group $^2E_6(2)$.
\end{abstract}

\section{Introduction}
\label{sec:intro}

Let $G$ be a group and $D$ a conjugacy class of involutions generating $G$ such that for all $d,e\in D$ 
the order of $de$ is equal to $1$, $2$ or $3$. Then $D$ is called a class  of {\em $3$-transpositions} in 
$G$. Groups generated by $3$-transpositions have been introduced by Fischer \cite{fischer} and studied 
by various authors since; e.g., see \cite{asch,vanbon,cuha,we1,we2}.

Given a class $D$ of $3$-transpositions in a group $G$, we define the {\em Fischer space} $\Pi(D)$
to be the partial linear space with $D$ as point set  and as lines the triples of points of the form 
$\{d,e,e^d=d^e\}$, where $d,e\in D$ are non-commuting. Thus, three $3$-transpositions on a line 
generate a subgroup  isomorphic to $Sym_3$, and vice versa, every subgroup $Sym_3$ containing 
involutions from $D$ produces a line.

The involutions from $D$ on two intersecting lines in a Fischer space generate a subgroup isomorphic to 
$\Sym_4$ or to a central quotient of the group $3^{1+2}\split 2$; e.g. see \cite{cuha}. The subspace
of the Fischer space generated by these two lines is then isomorphic to the dual of the affine plane of 
order $2$ or to the affine plane of order $3$, respectively. It was already noticed by Buekenhout, 
that  Fischer spaces are characterized by the property that any two intersecting lines generate such 
subspaces; see for example \cite{cuha}.

By $\F_2D$ we denote the $\F_2$ vector space on the set  of {\em finite} subsets of $D$, where addition 
of two sets is defined by the symmetric difference. We identify a point  $d\in D$ with the vector $\{d\}$ 
of this vector space. Note that this makes $D$ a basis of $\F_2D$, so we can write vectors from 
$\F_2D$ simply as linear combinations of the 3-transpositions from $D$.  Let $\A(D)$ be the algebra on 
$\F_2D$ whose multiplication $*$ is the linear expansion of the multiplication defined for $d,e\in D$ by:
\[d*e:=\begin{cases}
d+e+f \ \text{if} \ \{d,e,f\} \ \text{ is a  line}\\
0 \ \text{otherwise.}
\end{cases}
\]
In other words, if $d$ and $e$ commute (in particular, when $d=e$) then $d*e=0$ and if $d$ and $e$ do 
not commute then $d*e$ is the line that passes through them. The group $G=\langle D\rangle$ acts on 
$D$ by conjugation and it preserves lines. Hence it  induces  a group of automorphisms of both the Fischer 
space $\Pi(D)$ and the algebra $\A(D)$. Note that the action of $G$ is faithful only if $Z(G)=1$. 

We refer to the algebras $\A(D)$ and their $G$-invariant quotients as {\em $3$-transposition algebras}. 
Notice that Simon Norton  \cite{norton} considered a similar class of algebras, but defined over the reals.
Another related class of algebras (related to both Norton's algebras and our 3-transposition algebras) is
the class of Majorana algebras introduced by Alexander Ivanov \cite{ivanov}. In the (2B,3C)-case the Majorana
algebras are related to a subclass of groups generated by 3-transpositions. They, furthermore, admit a
natural basis with integral structure constants, and after reduction modulo two they produce our
3-transposition algebras.

The $3$-transposition algebra $\A(D)$ is endowed with a natural bilinear form defined as follows: for 
$d,e\in D$ we set $\form{d}{e}=0$, if $d$ and $e$ commute, and $\form{d}{e}=1$, if $d$ and $e$ do not 
commute. Since $D$ is a basis of $\A(D)$, this extends by linearity to the entire $\A(D)$. Note that 
$\form{d}{d}=0$ for every $d\in D$, that is, this bilinear form is symplectic. Another important property of 
the form $\form{\cdot}{\cdot}$ is that it {\em associates} with the algebra product, that is, 
\[\form{u}{v*w}=\form{u*v}{w} \]
for all $u,v,w\in\A(D)$. This will be verified in Proposition \ref{prop:associates}.

Clearly, the action of $G$ leaves the form invariant. Hence the radical $\V(D)$ of $\form{\cdot}{\cdot}$ is 
$G$-invariant. More interestingly, because the form associates with the algebra product, $\V(D)$ is an 
ideal of $\A(D)$! We call $\V(D)$ the {\em vanishing ideal} of $\A(D)$ and we call the elements of $\V(D)$, 
viewed as subsets of $D$, the {\em vanishing sets}.

Note that with the exception of the trivial case, where $G$ is a group of order two, the form 
$\form{\cdot}{\cdot}$ is nonzero, which means that $\V(D)$ is a proper ideal of $\A(D)$. As we are mainly 
interested in simple algebras, we will study the algebra $\OA:=\A(D)/\V(D)$, in which vanishing sets are 
reduced to zero. Note that in general $\OA$ does not need to be simple; however, it is very nearly so. We will 
show (see Proposition \ref{prop:invariant}) that every $G$-invariant proper ideal of $\A(D)$ is contained in 
$\V(D)$. In particular, $\OA$ is always semisimple.
 
The following relation $\tau$ plays a key role in the theory of Fischer spaces. For an element $d\in D$ we 
denote by $A_d$ the set of all points collinear to but distinct from $d$, i.e., $A_d:=\{e\in D\mid o(de)=3\}$.
If for $d,e\in D$ we have $A_d=A_e$, then we write $d\tau e$. The relation $\tau$ is an equivalence 
relation and it is related to the existence of a certain normal $2$-subgroup of $G$; see \cite{cuha}. Indeed, 
the  subgroup $\tau(G)=\langle de\mid d,e\in D, d\tau e\rangle$ is normal in $G$. The image $\bar D$ 
of $D$ in $\bar G=G/\tau(G)$ is  a class of $3$-transpositions of $\bar G$. Note that for any two elements 
$d,e\in D$ with $d\tau e$, we have $d+e\in\V(D)$ (that is, $\{d,e\}$ is a vanishing set in $D$). It follows that 
the algebras $\A(D)/\V(D)$ and $\A(\bar D)/\V(\bar D)$ are isomorphic. Thus, we can restrict our attention to 
the case where the relation $\tau$ is trivial. This is a strong condition that significantly simplifies the possible 
structure of $G$.
 
We will mainly focus on the case where the algebra $\OA$ is a Lie algebra. It
turns out that $\OA$ is a Lie algebra if and only if the affine planes of the
Fischer space are vanishing sets; see \PropRef{prop:islie}. This observation
makes it possible to determine which simple Lie algebras are (quotients of)
$3$-transposition algebras.

The main result of this paper is as follows. We use the Atlas notation for
groups; in particular, $p^n$ stands for the elementary abelian group of order
$p^n$. The colon $:$ indicates a split extension (semidirect product), where the
left side is normal. Furthermore, in all such extensions the natural action of
the complement on the normal subgroup is assumed. Finally, for a Dynkin diagram
$X_n$, $W(X_n)$ stands for the corresponding Weyl group.  

\begin{thm}\label{thm:main}
Let $G$ be a nonabelian group generated by a class $D$ of $3$-transpositions. Suppose the relation $\tau$ 
on $D$ is trivial and suppose further that $\OA$ is a simple Lie algebra.

Then, up to the center of $G$, we have one of the following:
\begin{enumerate}
\item $D$ is the unique class of $3$-transpositions in $G=3^n\split W(A_n)$; the algebra $\OA$ is isomorphic 
to the simple Lie algebra of type  $^2A_n(2)$.
\item $D$ is the unique class of $3$-transpositions in $G=3^n\split W(D_n)$; the algebra $\OA$ is isomorphic to 
the simple Lie algebra of type $^2D_n(2)$ for odd $n$ and of type $D_n(2)$ for even $n$.
\item $D$ is the unique class of $3$-transpositions in the group $G=3^n\split W(E_n)$ with $n\in \{6,7,8\}$; the 
algebra $\OA$ is isomorphic to the Lie algebra of type $^2E_6(2)$ ($n=6$), $E_7(2)$ ($n=7$) or $E_8(2)$ 
($n=8$).
\item $D$ is the  class of transvections in $G=\SU_{n+1}(2)$; the algebra $\OA$ is isomorphic to the simple Lie 
algebra of type  $^2A_n(2)$.
\item $D$ is one of the two classes of $126$ reflections of ${\rm O}^-_6(3)$ and $\OA$ is isomorphic to the 
simple Lie algebra of type $^2A_5(2)$.
\item $D$ is one of the two classes of $117$ reflections of ${\rm O}^+_6(3)$ and $\OA$ is isomorphic to the 
simple Lie algebra of type $D_4(2)$.
\item $D$ is the unique class of $360$ $3$-transpositions in $G=\POmega_8^+(2)\split\Sym_3$; the algebra
$\OA$ is isomorphic to the simple Lie algebra of type $D_4(2)$.
\item $D$ is the class of $351$ reflections of $+$-type in $G=^+\!\!\Omega_7^+(3)$ and $\OA$ is isomorphic to 
the simple Lie algebra of type  $^2E_6(2)$.
\item $D$ is the unique class of $3510$ $3$-transpositions in $\Fi{22}$ and $\OA$ is isomorphic to the simple Lie
 algebra of type $^2E_6(2)$.
\end{enumerate}
\end{thm}

For the notation used in part (h) we refer the reader to  Section \ref{sec:affine-vanish}.

The above result and its proof provide a geometric argument for the embedding of a central extension of 
$\POmega_6^-(3)$ into $\SU_6(2)$ and of $\Fi{22}$ into $^2E_6(2)$. Indeed, we obtain the following result.
 
\begin{cor}\label{cor:embed}
\begin{enumerate}
\item The group ${\PSU}_6(2)$ contains a subgroup isomorphic to $\POmega^-_6(3)$ generated by root
elements (i.e., elations).
\item The group $^2E_6(2)$ contains subgroups isomorphic to ${^+\Omega^+}_7(3)$ and to $\Fi{22}$ 
generated by root elements.
\end{enumerate}
\end{cor}

The embedding of $\Fi{22}$ into $^2E_6(2)$ was first established by Bernd Fischer. It led him to the discovery of 
the Baby Monster sporadic simple group and hence also to the discovery of the Monster, see \cite{ivanov}.

We also notice that with a bit of extra effort we could have included the case that the group $G$  is infinite.
Indeed, in \cite{cuha} all Fischer spaces, finite and infinite, have been classified. The infinite ones turn out to be 
limits of finite ones. So we only find the infinite dimensional,  finitary versions of the unitary and orthogonal Lie 
algebras as in (a), (b), and (d) of the conclusion of Theorem \ref{thm:main}, if we allow infinite groups $G$.

The organization of the paper is as follows. In Section \ref{sec:3-trans} we develop some general theory for 
$3$-transposition algebras. In particular, we prove the algebra $\A$ to be a Lie algebra if and only if every affine 
plane of the Fischer space is vanishing. In Section \ref{sec:cartan-normalizer} we prove various forms of simple 
Lie algebras of classical type $A$, $D$ or $E$ over the field $\F_2$ to be $3$-transposition algebras. In Section 
\ref{sec:affine-vanish} we start with the proof of our main result, Theorem \ref{thm:main}. We determine those
Fischer spaces that give rise to Lie algebras. In Section \ref{sec:proof-mainthm} we complete the proof of  
\ThmRef{thm:main}. In Section \ref{sec:comp-results} we present some additional computational results on the 
dimensions of (arbitrary) $3$-transposition algebras for relatively small groups. 

\bigskip\noindent
{\bf Acknowledgement.} The authors want to thank Ralf Gramlich for many fruitful discussions on the topic of this 
paper. The first and fourth authors also gratefully acknowledge the EPSRC grant that made possible the visit of
Hans Cuypers to the University of Birmingham and joint work on a number of projects, including this one.

\section{\boldmath $3$-transposition algebras}
\label{sec:3-trans}

Let $G$ be a group and $D$ a conjugacy class of $3$-transpositions generating $G$. By $\A$ we denote the 
corresponding $3$-transposition algebra $\A(D)$. We will study $\A$ and the action of $G$ on $\A$. Since the 
center $Z(G)$ of $G$ acts trivially on $D$ and since $de\not\in Z(G)$ for all $d,e\in D$, $d\neq e$, we can assume $Z(G)=1$ 
whenever convenient. 

Recall that the product in $\A$ is defined by 
\[d*e:=\begin{cases}
d+e+f \ \text{if} \ \{d,e,f\} \ \text{ is a  line}\\
0 \ \text{otherwise.}
\end{cases}
\]
and the symplectic form on $\A$ is defined by
\[\form{d}{e}:=\begin{cases}
1\ \text{if $d$ and $e$ are collinear and distinct,}\\
0\ \text{otherwise.}
\end{cases}
\]
We first check that $\form{\cdot}{\cdot}$ associates with the product. 

\begin{prop} \label{prop:associates}
For $v,u,w\in\A$ we have 
\[ \form{u}{v*w}=\form{u*v}{w}. \]
\end{prop}

\begin{proof}
It suffices to prove that $\form{u}{v*w}=\form{u*v}{w}$ for $u,v,w\in D$.

Suppose $u,v,w\in D$. If $u$ commutes with $v$, then $u*v=0$ hence $\form{u*v}{w}=0$. Thus we have to show 
that $\form{u}{v*w}=0$. This is clear if $w$ commutes with $v$, as then $w*v=v*w=0$. Otherwise, let 
$t=w^v=v^w$ be the third point on the line through $v$ and $w$. Note that $u$ commutes either with both $w$ 
and $t$, or with neither. In either case, we obtain $\form{u}{v*w}=0$. (Recall that $\A$ is defined over $\F_2$ 
and hence the values of $\form{\cdot}{\cdot}$ are in $\F_2$, too.)

By symmetry, we can now restrict our attention to the case where both $u$ and $w$ do not commute with $v$. If all 
three points are together on a line, then clearly $\form{u}{v*w}=0=\form{u*v}{w}$. So we can assume that the 
three points span a plane. If this plane is an affine plane of order 3, then $w$ is collinear to all three points on the line 
$u*v$ and $u$ is collinear to all three points on $v*w$. Hence $\form{u}{v*w}=1=\form{u*v}{w}$. If the plane is the 
dual affine plane of order 2, then $w$ is collinear to two points on the line $u*v$ and $u$ to two points on $v*w$. 
Therefore $\form{u}{v*w}=0=\form{u*v}{w}$.
\end{proof}

Let $\V=\V(D)$ be the radical of $\form{\cdot}{\cdot}$. 

\begin{cor} \label{cor:ideal}
The radical $\V$ is an ideal of $\A$.
\end{cor}

\begin{proof}
Since $\V$ is a linear subspace, it suffices to show that $\V$ is closed with respect to multiplication with elements of $\A$. Let $u\in\V$ and $v\in\A$. 
It follows from Proposition \ref{prop:associates} that for every $w\in\A$ we have that 
$\form{u*v}{w}=\form{u}{v*w}=0$, since $u$ is in the radical $\V$. Thus, $\form{u*v}{w}=0$ for all $w\in\A$, 
that is, $u*v$ is in $\V$.
\end{proof}

Recall that we call $\V$ the vanishing ideal and the elements of $\V$, viewed as subsets of $D$, vanishing sets. Clearly, 
a finite subset $X$ of $D$ is vanishing if and only if $X$ is perpendicular to every $d\in D$, that is, if for every $d\in D$, the number of elements of $X$ not commuting with $d$ is even. Every finite Fischer space does contain nonempty 
vanishing subsets. Indeed, if $D$ is finite, then $D$  itself is a vanishing set, as every point $d$ does not commute 
with an even number of points, two on every line through $d$. For each of the various types of Fischer spaces we 
can find more types of vanishing sets.

 We now concentrate on $G$-invariant ideals and quotients of $\A$. It is clear that the symplectic form 
$\form{\cdot}{\cdot}$ is preserved by the action of $G$ obtained by linearly extending the conjugation action of $G$ on $D$.
Thus $\V$ is $G$-invariant. We claim that 
$\V$ is the unique maximal among all $G$-invariant proper ideals of $\A$.

The following observation will be useful.

\begin{lemma}
Let $X$ be  a finite subset of $D$ and $d\in D$. Then 
\[ d*X=\form{d}{X}d+X+X^d.\]
\end{lemma}

\begin{proof}
The expression in the right side is linear in $X$, so one only needs to verify it for $X$ of size 1, in which case the claim 
follows directly from the definition of the product on $\A$.
\end{proof}

We now turn to the main claim.

\begin{prop}\label{prop:invariant}
Every $G$-invariant proper ideal of $\A$ is contained in $\V$.
\end{prop}

\begin{proof}
Suppose $\I$ is a $G$-invariant ideal containing a nonzero element $X$ not in $\V$. Then we can find an element $d\in D$ with $\form{d}{X}\neq 0$. Since $\I$ is an ideal, 
$d*X\in\I$. On the other hand, by the above lemma, $d*X=d+X+X^d$. Note that $X\in\I$ and also $X^d\in\I$, since
$\I$ is $G$-invariant. It follows that $d\in\I$. However, now by $G$-invariance, $d^G=D$ is contained in $\I$ and so 
$\I=\A$.
\end{proof}

The above result has the following important consequence. Recall that $\OA=\A(D)/\V(D)$.

\begin{prop} \label{prop:semisimple}
If $D$ is finite then $\OA$ is the direct product of isomorphic simple algebras. In particular,
$\OA$ is semisimple.
\end{prop}

\begin{proof}
If $D$ is finite then $\A$ is finite dimensional and hence it contains a minimal ideal $\I$. For
$g\in G$, if $\I^g\neq\I$ then $\I\cap\I^g=0$ and hence $\I\I^g=0$. It follows from here 
inductively that the orbit of $\I$ under $G$ generates an ideal $\J$ that is a direct sum of 
several conjugates of $\I$. Since $\J$ is manifestly $G$-invariant, Proposition 
\ref{prop:invariant} implies that $\J=\OA$. Finally, since $\I$ is a minimal ideal and since 
every factor of $\OA$, other than $\I$ itself, annihilates $\I$, it follows that $\I$ is a simple 
algebra.
\end{proof}

Note that we were careful in this proof not to state that $\OA$ is the direct product of all 
ideals $\I^g$. This is, however, almost always the case. Indeed, if some conjugate $\I^g$ 
does not appear in the direct product decomposition then every factor annihilates 
$\I^g$, which clearly means that $\I^g$ (and hence also $\I$ and the entire $\OA$) is a 
trivial algebra. It easily follows from the definition of the product of $\A$ that $\OA$ is 
trivial if and only if every line in the Fischer space $\Pi(D)$ is vanishing.  An example of a Fischer
space satisfying this property is the dual affine plane of order 2. For more examples 
see Lemma \ref{lem:symplectic}.

The above discussion yields the following.

\begin{cor}
If the finite Fischer space $\Pi(D)$ has at least one line that is not vanishing then $\OA$ 
is the direct product of all its minimal ideals. Furthermore, $G$ transitively permutes
the minimal ideals of $\OA$.
\end{cor}

We also note the following important case.

\begin{prop}\label{prop:simple}
Suppose $A$ is a simple quotient algebra of $\A$. If $G$ induces a group of automorphisms on $A$,
then $A$ is isomorphic to $\OA$. 
\end{prop}

\begin{proof}
Clearly, $A$ is isomorphic to $\A/\I$ for some ideal $\I$. Since $G$ induces an  action on $A$, the ideal 
$\I$ is $G$-invariant. However, in this case, by Proposition \ref{prop:invariant}, $\I$ is contained in $\V$. 
Simplicity of $A$ implies now that $\I=\V$.
\end{proof}

The following lemma will be used in Section \ref{sec:proof-mainthm}.

\begin{lemma}\label{lem:subalgebra}
Let $E$ be a subspace 
of the Fischer space $\Pi(D)$ (i.e., $E$ satisfies $d^e\in E$ for all $d,e\in E$). Then the algebra $\OA$   
contains a subalgebra that has a quotient isomorphic to $\A(E)/\V(E)$.
\end{lemma}

\begin{proof}
The elements of $E$ generate a subalgebra $\A(E)$ of $\A(D)$. As every vanishing subset of $D$ which 
is contained in $E$ is also a vanishing subset of $E$, we find that $\V(D)\cap\A(E)\subseteq \V(E)$.
Hence, $(\A(E)+\V(D))/\V(D)$ is the subalgebra we are looking for.
\end{proof}

Finally, we need a criterion which allows us to decide when $\OA$ is a Lie algebra.

\begin{prop}\label{prop:islie}
Let $\I$ be an ideal of $\A$. Then $\A/\I$ is a Lie algebra, if and only if every plane of $\Pi(D)$ isomorphic 
to the affine plane of order $3$ is in $\I$.
\end{prop}

\begin{proof}
Clearly the Jacobi identity holds in $\A$ if and only if it holds for any three elements of $D$.
Let $d,e,f\in D$. We will check the Jacobi identity for these elements,  that is, the equality
\[ (d*e)*f+(e*f)*d+(f*d)*e=0. \]

If two of the three elements are equal, say $d=e$, then 
\[\begin{array}{lll}
(d*e)*f+(e*f)*d+(f*d)*e&=&(d*d)*f+(d*f)*d+(f*d)*d\\
                        &=& 0+(d*f)*d+(d*f)*f\\
                        &=& 0.\\
\end{array}\]
Hence, we can assume that the three elements are pairwise distinct.

If $\{d,e,f\}$ is a line, then $(d*e)*f=(d+e+f)*f=d+e+f+d+e+f=0$. By symmetry,
$(d*e)*f+(e*f)*d+(f*d)*e=0+0+0=0$. Thus, assume that the elements $d,e,f$ are not 
on a single line.

If there is a point (say $d$) among $d,e,f$ which is neither collinear to $e$ nor to $f$ 
(i.e., $d$ commutes with both $e$ and $f$), then all three terms $(d*e)*f$, $(e*f)*d$ and 
$(f*d)*e$ are $0$. Indeed, not only $d*e=f*d=0$, but also $(e*f)*d=0$, since either
$e*f=0$ or else $e*f=e+f+t$ for $t=e^f=f^e$ and so $(e*f)*d=e*d+f*d+t*d=0+0+0$, 
as $d$ is not collinear to $t$ (indeed, if $d$ commutes with $e$ and $f$ then it also 
commutes with $t=e^f$). So the  Jacobi identity holds also in this case.

This leaves the situation where $d,e,f$ are three points in a plane of the Fischer space 
$\Pi(D)$. Moreover, we can assume that $e$ is collinear to both $d$ and $f$. First suppose 
that $d,e,f$ are in a dual affine plane. If $d*f=0$, then $(e*f)*d$ equals the sum of the 
two lines on $d$. Also, $(f*d)*e=0$ and $(d*e)*f$ equals the sum of the two lines on $f$.
However, the sum of the two lines on $d$ and the two lines on $f$ equals zero, which 
establishes the Jacobi identity in this case.

If $d*f\neq 0$ then $(d*e)*f+(e*f)*d+(f*d)*e$ is the sum of the lines of the plane on $d$, 
$e$ and $f$. Again this sum is equal to zero.

Now assume that $d,e,f$ are three points inside an affine plane, say $\pi$. Then $(d*e)*f$ 
is equal to the sum of the three lines in the plane on $f$ meeting the line through $d$ and 
$e$. But then $(d*e)*f+(e*f)*d+(f*d)*e$, after canceling the lines through $d$ and $e$, 
through $d$ and $f$, and through $e$ and $f$, each of which appears in the sum twice, 
equals the sum of the three lines passing through $f$ and $d^e$, through $d$ and $e^f$,
and through $e$ and $f^d$. This is the sum of three parallel lines in the plane $\pi$ and 
hence it equals to $\pi$ itself.

Thus, the Jacobi identity holds in $\A/\I$, if and only if every such $\pi$ is in $\I$. 
\end{proof}

Clearly, this gives us the following.

\begin{cor} \label{cor:criterion}
$\OA$ is a Lie algebra if and only if every affine plane in $\Pi(D)$ is a vanishing subset 
of $D$.
\end{cor}

\section{\boldmath The normalizer of a Cartan subalgebra as a $3$-transposition group}
\label{sec:cartan-normalizer}

In this section we describe some examples of classical Lie algebras obtained
as quotients of $3$-transposition algebras.

Let $\mathfrak{g}$ be a split classical Lie algebra over $\F_4$, the field
with $4$ elements, of simply laced type, i.e., of type $A_n$ ($n\geq
1$), $D_n$ ($n\geq 4$) or $E_n$ ($n=6,7$, or $8$). Let $\mathfrak{h}$ be a
Cartan subalgebra of $\mathfrak{g}$  and $\Phi$ the corresponding root
system. We can then decompose $\mathfrak{g}$ as the sum of the Cartan
subalgebra $\mathfrak{h}$ and the corresponding root spaces:
$$\mathfrak{g}=\mathfrak{h}\oplus\bigoplus_{\alpha\in\Phi}\mathfrak{X}_\alpha,$$
where $\mathfrak{X}_\alpha$ denotes the root space associated to the root
$\alpha\in\Phi$.

For each $\alpha\in\Phi$ let $\mathfrak{x}_\alpha$ be a nonzero element in
$\mathfrak{X}_\alpha$, such that the elements $\mathfrak{x}_\alpha$ together
with the root elements
$\mathfrak{h}_\alpha=[\mathfrak{x}_\alpha,\mathfrak{x}_{-\alpha}]\in\mathfrak{h}$ form
a Chevalley basis for $\mathfrak{g}$.

The normalizer  of $\mathfrak{h}$ in the automorphism group of $\mathfrak{g}$ contains a
subgroup $G$ which is a split extension of an elementary abelian group $H\simeq
(\F_4^*)^n$ of order $3^n$ by a group $W$ isomorphic to the Weyl group of
$\mathfrak{g}$. This group $G$ is generated by its conjugacy class $D$ containing
the class of reflections from $W$. Moreover, this class of involutions is a
class of $3$-transpositions in $G$; see also \cite{cuha}, examples {\bf PR 2},
{\bf PR 9}, {\bf PR 10} and {\bf PR 11}. We now describe a Fischer space on a
set of elements of $\mathfrak{g}$ isomorphic to the Fischer space of $D$.

Let $\Delta$ be the set of elements of the form
$$x(\alpha,\omega):=\omega\mathfrak{x_\alpha}+\overline{\omega}\mathfrak{x}_{-\alpha}+\omega\overline{\omega}[\mathfrak{x}_\alpha,\mathfrak{x}_{-\alpha}],$$
where $\alpha\in\Phi$, and $\omega$ is an element of $\F_4$ with
$\overline\omega=\omega^2$ being the image of $\omega$ under the unique
nontrivial field automorphism of $\F_4$.

The elements in $\Delta$ are all long root elements in $\mathfrak{g}$, as one
easily checks within $\langle\mathfrak{X}_\alpha,\mathfrak{X}_{-\alpha}\rangle$.
(Actually, they form a $G$-orbit of root elements.) Moreover, they linearly span
$\mathfrak{g}$. As the bar automorphism of $\F_4$ acts on $\Delta$ as a
permutation, mapping $x(\alpha,\omega)$ to $x(\alpha,\overline{\omega})$, the
elements of $\Delta$ also span the Lie algebra $\mathfrak{g}_{2}$ over $\F_2$
that they generate. Notice that this subalgebra is also fixed by the Chevalley
involution of $\mathfrak{g}$. The  algebra  $\mathfrak{g}_{2}$ is, up to its
center, a simple Lie algebra; its type is given in the table below.

\bigskip
\begin{tabular}{|l|l|l|}
\hline
Type $\mathfrak{g}$ &Type $\mathfrak{g}_2$ & dimension $\mathfrak{g}_2/Z(\mathfrak{g}_2)$\\
\hline
$A_n$          & $^2A_n$ & $(n+1)^2-2$ if $n$ is odd, $(n+1)^2-1$ if $n$ even\\ 
$D_n$ ($n$ odd)& $^2D_n$&   $2n^2-n-1$ \\ 
$D_n$ ($n$ even) & $D_n$  & $2n^2-n-2$ \\
$E_6$          & $^2E_6$&  $78$\\ 
$E_7$          & $E_7$  &  $132$ \\ 
$E_8$          & $E_8$  &  $248$\\ 
\hline
\end{tabular}

\medskip
We claim that the elements of $\Delta$ form a Fischer space isomorphic
to the Fischer space $\Pi(D)$ in $G$. The lines of this Fischer space
are the triples $\{\mathfrak{x},\mathfrak{y},\mathfrak{z}\}$ of elements from
$\Delta$ satisfying $[\mathfrak{x},\mathfrak{y}]=\mathfrak{x}+\mathfrak{y}+\mathfrak{z}$.

Let $\alpha$ and $\beta$ be two distinct roots in $\Phi$ and
$\mu,\nu\in\F_4^*$. We want to determine 
$$[x(\alpha,\mu),x(\beta,\nu)].$$
Since $x(\beta,\nu)=x(-\beta,\overline{\nu})$ we only have to consider
the following three cases: $\alpha=\beta$, or $\alpha+\beta\in\Phi$, or
both $\alpha+\beta$ and $\alpha-\beta$ are nonzero and not in $\Phi$.

If $\alpha=\beta$, then we have
\[
    [x(\alpha,\mu),x(\beta,\nu)]
  = \mathfrak{h}_\alpha
  = x(\alpha,\mu)+x(\alpha,\nu)+x(\alpha,\mu+\nu)
  .
\]

If $\alpha+\beta\in\Phi$, then 
\[
  [x(\alpha,\mu),x(\beta,\nu)]
= x(\alpha,\mu)+x(\beta,\nu)+x(\alpha+\beta,\mu\nu).
\]

Finally, if both $\alpha+\beta$ and $\alpha-\beta$ are not in $\Phi$, then
$$[x(\alpha,\mu),x(\beta,\nu)]=0.$$

This implies that $\Delta$ is closed under the lines defined above. In
particular, these lines define a partial linear space on $\Delta$. We prove this
space to be a Fischer space. Consider the subgeometry spanned by two
intersecting lines. Up to replacing a root by its negative, this leads to the
following two cases. First, the two lines are spanned by the points
$x(\alpha,\mu)$ and $x(\beta,\nu)$, and $x(\gamma,\rho)$ and $x(\beta,\nu)$,
respectively, where $\alpha$, $\beta$ and $\gamma$ are three distinct positive
roots, with $\alpha$ and $\gamma$ orthogonal, and $\alpha+\beta$, $\beta+\gamma$
also positive roots, and $\mu,\nu,\rho\in\F_4^*$. In this case
$\alpha+\beta+\gamma$ is also a root, and we find the following six elements of
$\Delta$ to form a dual affine plane of order 2 containing the two lines:
$$x(\alpha,\mu), x(\beta,\nu), x(\alpha+\beta,\mu\nu),$$
$$x(\gamma,\rho), x(\beta+\gamma,\nu\rho), x(\alpha+\beta+\gamma,\mu\nu\rho).$$

In the second case, the two lines are spanned by the points $x(\alpha,\mu)$ and
$x(\beta,\nu)$, and $x(\alpha,\rho)$ and $x(\beta,\nu)$, respectively, where
$\alpha$ and $\beta$ are distinct positive roots, with $\alpha+\beta$ also
being a positive root.

In this case the nine points of the form $x(\delta,\tau)$, where
$\delta\in\{\alpha,\beta,\alpha+\beta\}$ and $\tau\in\F_4^*$, form a subspace
of $\Delta$ isomorphic to an affine plane of order 3.

So, indeed, $\Delta$ carries the structure of a Fischer space, which, as is
clear from the description of $\Delta$, is isomorphic to the Fischer space on
the set $D$ of $3$-transpositions in $W$.

This implies that the algebra $\mathfrak{g}_{2}$ is a quotient of the
$3$-transposition algebra $\A(D)$. So, with the above notation we have proven
the following result.

\begin{prop}\label{prop:weyl}
Suppose $D$ is the class of $3$-transpositions in the group $G$, as defined in
this section. Then the Lie algebra $\A(D)/\V(D)$ is isomorphic to the Lie
algebra $\mathfrak{g}_2$ modulo its center.
\end{prop}

\begin{proof}
By construction, the group $G$ is contained in the automorphism group of
$\mathfrak{g}_2$. Since the latter Lie algebra modulo its center is simple, we can
apply \PropRef{prop:simple}.
\end{proof}

\section{Fischer spaces in which affine planes vanish}
\label{sec:affine-vanish}

The $3$-transposition algebras giving rise to Lie algebras come from Fischer
spaces in which affine planes are vanishing; see \PropRef{prop:islie}. Thus, in
order to find all Lie algebras among the $3$-transposition algebras, we can
restrict our attention to Fischer spaces in which affine planes are vanishing.
It is the purpose of this section to classify all such spaces. For this we make
use of (parts of) the classification of Fischer spaces as presented in
\cite{cuha}.

We introduce some notation. Denote by $D_d$ the set of all points in $D$ not
collinear with $d$, i.e., $D_d:=\{e\in D\mid o(de)=2\}$ (observe that $D=A_d\cup
D_d\cup\{d\}$). As we already saw in the introduction, the relation $\tau$ on
$D$ defined by $d\tau e$ for $d,e \in D$ if and only if $A_d=A_e$ is related to
the existence of a normal $2$-subgroup in $G$. The relation $\theta$ on $D$,
defined by $d\theta e$ if and only if $D_d=D_e$ is related to normal
$3$-subgroups. Indeed, the subgroup $\theta(G)=\langle de\mid d,e\in D,d\theta
e\rangle$ is normal in $G$ and it is a 3-group.

The space $\Pi$ is said to be of {\em symplectic} type if it contains a dual
affine plane, but no affine planes. It is called of {\em orthogonal} type if it
does contain an affine plane, but every point not in the plane is collinear with
0, 6 or all points in the plane. This excludes the case where $G$ contains a
subgroup generated by elements from $D$ isomorphic to a central quotient of
$2^{1+6}\split\SU_3(2)$. Such subgroups appear in the unitary groups over
$\F_4$. Thus, continuing in this vein, we say that $\Pi$ is of {\em unitary}
type when it contains an affine plane and a point outside of the plane that is
collinear to exactly $8$ points of the plane, but does not contain a subspace
isomorphic to the Fischer space of the unique class of $3$-transpositions in
$\POmega^+_8(2)\split\Sym_3$ (a subgroup of all five sporadic examples from
\cite{cuha}). Finally if $\Pi$ does  contains a subspace isomorphic to the
Fischer space of $\POmega^+_8 (2)\split\Sym_3$, then it is of {\em sporadic}
type.

The Fischer space $\Pi(D)$ (as well as the group $G=\langle D\rangle$) will be
called {\em irreducible} if and only if $\Pi(D)$ is connected, and both
relations $\tau$ and $\theta$ are trivial. We recall the first main theorem
from \cite{cuha}.

\begin{thm}[{\cite[Theorem 1.1]{cuha}}]\label{thm:classthm}
Let $G$ be a group generated by a conjugacy class $D$ of $3$-transpositions.
If $G$ is irreducible, then, up to a center, we may identify $D$ with one of
the following:
\begin{enumerate} 
\item
the transposition class of a symmetric group;

\item
the transvection class of the isometry group of a nondegenerate
orthogonal space over $\F_{2}$;

\item
the transvection class of the isometry group of a nondegenerate
symplectic space over $\F_{2}$;

\item
a reflection class of the isometry group of a nondegenerate
orthogonal space over $\F_{3}$;

\item
the transvection class of the isometry group of a nondegenerate
unitary space over $\F_{4}$;

\item
a unique class of involutions in one of the five groups
$\POmega_8^+(2)\split\Sym_3$, $\POmega_8^+(3)\split\Sym_3$, $\Fi{22}$,
$\Fi{23}$, or $\Fi{24}$.
\end{enumerate}
\end{thm}

The $3$-transposition classes of cases (a)--(c) of \ThmRef{thm:classthm} are of
symplectic type, those described in case (d) of orthogonal type, the ones
described in (e) of unitary type, and finally, those in (f) of sporadic type.

In all cases except for (d), the class of $3$-transpositions is unique. An
orthogonal group over $\F_3$, however, contains two classes of reflections. Let
$(V,Q)$ be a nondegenerate orthogonal space of dimension $n$ over $\F_3$ and
suppose $f$ is the associated bilinear form. If $n$ is finite, then, up to
isometry, there are two choices for the form $Q$, distinguished by their
discriminant $\delta=\pm 1$. This discriminant is, in even dimension, determined
by the Witt sign $\epsilon$ of $Q$, which is defined as $+1$ if the Witt index
(that is, the dimension of maximal isotropic subspaces) equals $\frac{n}{2}$ and
as $-1$ if the Witt index equals $\frac{n}{2}-1$. Indeed, for even $n$ we have
$\epsilon\delta=(-1)^{\frac{(n+1)n}{2}}$. We use this formula to define the Witt
sign $\epsilon$ also in odd dimensions. We write $\epsilon=\pm$ for
$\epsilon=\pm 1$. For all $n$, we denote $\O(V,Q)$ by $\O^\epsilon_n(3)$ and its
derived subgroup by $\Omega^\epsilon_n(3)$. For odd $n$, the $\epsilon$ is often
left out, since then $\O^+_n(3)$ and $\O^-_n(3)$ are isomorphic.

For each vector $x\in V$ with $Q(x)\neq 0$, the reflection $r_x:v\mapsto
v+f(v,x)Q(x)x$ is an element of the orthogonal group $\O(V,Q)$. There are two
classes of reflections in $\O(V,Q)$; those of $+$-type with $Q(x)=1$ and those
of minus type with $Q(x)=-1$.

By $^\gamma\Omega^\epsilon_n(3)$ we denote the subgroup of $\O(V,Q)$ generated
by the reflections of type $\gamma$.

The aim of this section is to prove the following result.

\begin{prop}\label{prop:class}
Let $\Pi$ be a finite connected Fischer space with trivial relation $\tau$. Then 
the affine planes in $\Pi$ are vanishing, if and only if one of the following 
holds:
\begin{enumerate}
\item $\Pi$ is of symplectic type.
\item $\Pi$ is isomorphic to the Fischer space on the unique class of
$3$-transpositions in $3^{n}\split W(X_n)$, where $X_n$ is $A_n$ $(n\geq 2)$ ,
$D_n$ $(n\geq 4)$ or $E_n$ $(n=6,7,8)$.
\item $\Pi$ is isomorphic to the Fischer space of transvections in the group
$\SU_n(2)$ with $n\geq 3$.
\item $\Pi$ is isomorphic to the Fischer space of the class of $126$ 
$+$-reflections in ${^+\Omega}^-_6(3)$.
\item $\Pi$ is isomorphic to the Fischer space of either one of the two classes
of $117$ $+$-reflections in ${^+\Omega}^+_6(3)$.
\item  $\Pi$ is isomorphic to the Fischer space on the  $351$  $+$-reflections
in ${^+\Omega^+}_7(3)$.
\item $\Pi$ is isomorphic to the Fischer space of the unique class of
$360$ $3$-transpositions of $\POmega_8(2)\split\Sym_3$.
\item $\Pi$ is isomorphic to the Fischer space of  the unique class of 
$3$-transpositions  of $\Fi{22}$.
\end{enumerate}
\end{prop}

The proof of this proposition is given in the remainder of this section.

Assume that $\Pi$ is a Fischer space of some class $D$ of $3$-transpositions in
a group $G=\langle D\rangle$. Moreover, assume $\Pi$ satisfies the hypothesis of
\PropRef{prop:class}.

If $\Pi$ is of symplectic type, there is nothing to prove. So, we assume that
$\Pi$ does contain affine planes. The following  observation is very useful.

\begin{lemma}\label{lem:no3space}
A Fischer space in which affine planes are vanishing does not contain affine 
$3$-spaces.
\end{lemma}

\begin{proof}
Inside an affine $3$-space, an affine plane is not vanishing.
\end{proof}

As already noticed in the introduction, we focus on the Fischer spaces for which
the relation $\tau$ is trivial. However, let us first assume that not only
$\tau$, but also $\theta$ is trivial.

Suppose $\Pi$ is of orthogonal type. In this case $D$ can be identified with a
class of reflections in an orthogonal group $\O(W,q)$ for some nondegenerate
orthogonal space $(W,q)$  over the field $\F_3$, as in case (d) of 
\ThmRef{thm:classthm}. However, since by \LemRef{lem:no3space} there are no
affine $3$-spaces in $\Pi$, the dimension of $W$ is restricted by $5\leq\dim W
\leq7$.

If $\dim(W)=5$, then only one class of reflections provides a Fischer space with
affine planes, namely those reflections, which centralize an orthogonal
$4$-space of maximal Witt index. This Fischer space is then isomorphic to the
Fischer space of the unique class of $3$-transpositions in $\SU_4(2)$. Hence we
are in case (c) of the proposition.

If $\dim(W)=6$, then we have case  (d) or (e) of \PropRef{prop:class}.

Finally, if $\dim(W)=7$, then by the condition that a point always centralizes
at least a line in a plane not containing the point, only the class of
reflections centralizing an orthogonal $6$-space with non-maximal Witt index
satisfies our assumptions. This class is the class of $+$-reflections in the
group $^+\Omega^+_7(3)$, as described under (f) of \PropRef{prop:class}.

Next assume that $\Pi$ is of unitary type. Then we are in case (e) of
\ThmRef{thm:classthm} and hence in case (c) of \PropRef{prop:class}. The fact
that the affine planes inside all these unitary spaces are vanishing can easily
be checked within $\SU_4(2)$ and $\SU_5(2)$.

It remains to consider the case where the Fischer space is sporadic, as in case
(g) of \ThmRef{thm:classthm}.

The commuting graph of the Fischer space of $\POmega_8(3)\split\Sym_3$ consists
of three parts, each forming a Fischer space for $\POmega_8(3)$. Points of one
part are collinear to all points of the other two parts. As each part contains 
affine planes, we conclude that these planes are not vanishing. Since the
Fischer spaces related to $\Fi{23}$ and $\Fi{24}$ contain subspaces isomorphic
to the Fischer space of $\POmega_8(3)\split\Sym_3$, they contain non-vanishing
affine planes, too.

The Fischer spaces of $\POmega_8(2)\split\Sym_3$ and $\Fi{22}$ occur in (f) and
(g) of \PropRef{prop:class}, respectively. It remains to prove that affine
planes in these Fischer spaces are vanishing.

\begin{lemma}
Let $\Pi$ be the Fischer space of $\POmega_8(2)\split \Sym_3$. Then all affine
planes in $\Pi$ are vanishing.
\end{lemma} 

\begin{proof}
The non-collinearity graph of $\Pi$ partitions into three parts, each part being
a Fischer space for $\POmega_8(2)$. The latter space is of symplectic type and
does not contain affine planes.

Therefore, every affine plane $\pi$ of $\Pi$ meets each part of the partition in
a line. The three lines form a parallel class of lines in the plane. If $d$ is a
point of $\Pi$, then $d$ centralizes one or three points of the line of $\pi$
that lies in the same part as $d$, and no other point of $\pi$. Hence $\pi$ is
vanishing. \end{proof}

Now assume $\Pi$ to be the Fischer space of $\Fi{22}$. To show that the affine
planes in $\Pi$ are vanishing, we consider its subspace on the $693$ points not
collinear to a fixed point $d$. This latter space is the Fischer space related
to $\SU_6(2)$; e.g., see \cite{asch,vanbon,cuha}. This space has the following
properties.

\begin{lemma}\label{lem:u6}
In the Fischer space of  $\SU_6(2)$ we have the following:
\begin{enumerate}
\item each point is on $256$ lines;
\item each line is in $40$ affine planes and $135$ dual affine planes;
\item there are in total $59136$ lines and $197120$ affine planes;
\item each point lies in $2560$ affine planes.
\end{enumerate}
\end{lemma}

\begin{proof}
The collinearity graph of the Fischer space of  $\SU_6(2)$ is strongly regular
with parameters $(v,k,\lambda,\mu)=(693,512,376,384)$ (see for example
\cite{asch,vanbon}). So, each point in this Fischer space is collinear to $512$
points. As any pair of collinear points determines a unique line, there are
$512/2=256$ lines per point. This proves (a).

Fix collinear points $d$ and $e$. Then each dual affine plane on the line $l$
through $d$ and $e$ contains a unique point collinear to $e$, but not to $d$.
Moreover, each point with this property determines a unique dual affine plane
containing $l$. Hence, there are $k-\lambda-1=135$ dual affine planes on $d$ and
$e$. Inside these planes we find $1+135=136$ common neighbors of $d$ and $e$.
So, the remaining $376-136=240$ common neighbors of $d$ and $e$ are in
$240/6=40$ affine planes, proving (b).

Double counting now implies (c) and (d).
\end{proof}

The collinearity graph of the Fischer space of  $\Fi{22}$ is strongly regular
with parameters $(v,k,\lambda,\mu)=(3510,2816,2248,2304)$; see for example
\cite{asch,vanbon} or \cite{cuha}. By similar arguments as used in the proof of
the previous lemma we obtain:

\begin{lemma}\label{lem:f22}
In the Fischer space of  $\Fi{22}$ we have the following:
\begin{enumerate}
\item each point is on $1408$ lines;
\item each line is in $280$ affine planes and $567$ dual affine planes;
\item there are in total $1647360$ lines and $38438400$ affine planes;
\item each point lies in $98560$ affine planes.
\end{enumerate}
\end{lemma}

We are now in a position to prove the following:

\begin{lemma}
The affine planes in the Fischer space of $\Fi{22}$ are vanishing.
\end{lemma}

\begin{proof}
Fix a point $d$. As we have seen above, there are $98560$ affine planes
containing $d$, and by \LemRef{lem:u6}, there are $197120$ affine planes
containing no point collinear to $d$.

Now fix a point $e$ not collinear to $d$ and a line $l$ through $e$ consisting
of points not collinear to $d$. By \LemRef{lem:u6} there are $693$ such points
$e$ and $256$ such lines on each $e$. Of the $280$ affine planes containing $l$,
there are $40$ planes having no point collinear to $d$. The other $240$ planes
contain $6$ points collinear to $d$. So, $e$ is in ${256\cdot 240}=61440$ affine
planes containing just a line of points not collinear to $d$. As $e$ lies in
$2560$ affine planes that contain no point collinear with $d$ (see
\LemRef{lem:u6}), there are $98560-61440-2560=34560$ planes through $e$
containing $8$ points collinear with $d$. So, we find $\frac{693\cdot 256\cdot
240}{3}=14192640$ affine planes containing just a line of points not collinear
to $d$, and $693\cdot 34560=23950080$ planes containing a unique point not
collinear to $d$.

However, by now we have accounted for $98560+197120+14192640+23950080=38438400$
affine planes. As this is in fact the total number of affine planes in $\Pi$, we
conclude that the point $d$ is collinear to $8$, $6$, or $0$ points of any
affine plane. This finishes the proof.
\end{proof}

We have covered the case where the relation $\theta$ on the set $D$ is trivial.
Now assume $\theta$ to be nontrivial. In this case the quotient $G/\theta(G)$ is
generated by a class $\overline{D}$ of $3$-transpositions whose Fischer space
does not contain affine planes. For otherwise, we would find affine $3$-spaces
in the space on $D$, contradicting \LemRef{lem:no3space}. Moreover, each
$\theta$-equivalence class consists of exactly $3$ elements, or else the Fischer
space on $\overline{D}$ contains no lines. In the latter case, $\Pi$ is an
affine plane, the Fischer space of $\SU_3(2)'$, as in case (c) of
\PropRef{prop:class}.

So we can assume that the Fischer space $\Pi(\overline{D})$ has lines but no
affine planes, and furthermore, that every $\theta$-equivalence class has size
3. This implies that, in the notation of \cite{cuha}, the Fischer space $\Pi(D)$
is  of orthogonal type with the property that in each affine plane $\pi$ there
is a line $l$ such that each point $d$ collinear with a point of $\pi$ is also
collinear to a point on $l$. Such spaces are classified in \cite[Theorem
6.13]{cuha}. By this theorem and the condition that each $\theta$-equivalence
class has size $3$ we find the examples described in case (b) of
\PropRef{prop:class}.

Indeed, if, in the notation of \cite[Theorem 6.13]{cuha}, $G$ is of type {\bf
PR1}, then the absence of an affine $3$-space leads to the case where $\Pi$ is
the Fischer space related to $3^n\split W(A_n)$ with $n=1,2$. If $G$ is of type
{\bf PR5}, {\bf PR9}, or {\bf PR10}, then we find the cases where $\Pi$ is the
Fischer space of $3^n\split W(E_n)$ with $n=6$, $7$, or $8$, respectively.
Hence, there only remains the case where $G$ is isomorphic to a group
$W(K,\Omega)$, a subgroup of the wreath product of a strong $\{2,3\}$-group $K$
and $\Sym_\Omega$, with $\Omega$ a set of size at least $4$, as described in
{\bf PR2} of \cite{cuha}. If $K$ is a $2$-group, then the corresponding Fischer
space does not contain affine planes. If $K$ contains subgroup of order $3^n$
with $n\geq 2$, then the Fischer space will contain affine $3$-spaces, which is
against our assumptions. Thus, $K$ is either cyclic of order $3$ or isomorphic
to $\Sym_3$. The first case leads to the Fischer spaces related to $3^n\split
W(A_n)$, $n\geq 4$, the second case to the Fischer spaces related to $3^n\split
W(D_n)$, $n\geq 4$, as described in case (b) of \PropRef{prop:class}.

\section{Proof of Theorem \ref{thm:main}}
\label{sec:proof-mainthm}

In this section we prove the main result of this paper, Theorem \ref{thm:main}.
We keep the notation as in the previous sections. In particular,
$\OA=\A(D)/\V(D)$.

\begin{lemma} \label{lem:symplectic}
If the Fischer space $\Pi(D)$ is symplectic, then $\OA$ is an abelian Lie
algebra. 
\end{lemma}

\begin{proof}
In this case the lines of the Fischer space are vanishing sets. Indeed, any
point $d\in D$ is collinear to $0$ or $2$ points (different from $d$) on any
line. So in $\OA$ the product of any two elements is $0$.
\end{proof}
 
The next lemma covers part (d) of \ThmRef{thm:main}.

\begin{lemma}
Suppose $D$ is the  class of transvections in the unitary group $\SU_{n+1}(2)$,
$n\geq 2$. Then $\OA$ is isomorphic to the Lie algebra of type $^2A_n(2)$ modulo
its center.
\end{lemma}

\begin{proof}
The long root elements in the Lie algebra of type $^2A_n(2)$ modulo its center
satisfy the same relations as the elements of $D$. Indeed, we can identify these
long root elements with the rank $1$ matrices in the unitary Lie algebra
$\mathfrak{su}_{n+1}(2)$, which are in a one-to-one correspondence with the
transvections in  $\SU_{n+1}(2)$. In view of simplicity of the Lie algebra of
type $^2A_n(2)$ modulo its center, the result follows from \PropRef{prop:simple}.
\end{proof}

We now turn our attention to the exceptional cases (e--i) of  \ThmRef{thm:main}.

\begin{lemma}
Suppose $D$ is one of the two classes of reflections in $\O_6^-(3)$. Then the
algebra $\OA$ is isomorphic to the Lie algebra of type $^2A_5(2)$ modulo its
center.
\end{lemma}

\begin{proof}
Let $(M,Q)$ be a $6$-dimensional orthogonal space over $\F_3$ of Witt index
$-1$. Up to isomorphism, we can identify the elements of $D$ with the
non-isotropic $1$-spaces $\langle m\rangle$ in $M$, with $Q(m)=1$. Every subset
of $D$ of pairwise commuting elements is contained in a set of size $6$
corresponding to an orthonormal basis of $M$. If $d\in D$ and  $B$ is such a set
of $6$ pairwise commuting elements corresponding to an orthonormal basis, then,
as simple computations in the orthogonal space reveal, $d$ is either contained
in $B$, or it commutes with exactly two elements of $B$. This implies that $B$
is a vanishing set. Moreover, it also implies that the set $\mathcal{C}$ of
subsets of 6 pairwise commuting elements of $D$ induces the structure of an
extended generalized quadrangle on $D$. In \cite{cush}, it has been shown that
the subspace of $\F_2D$ generated by the elements of $\mathcal{C}$ has
codimension $35$. As this subspace does not contain affine planes (which are of
odd weight), the vanishing ideal $\V$ of $\A$ has codimension at most $34$.

We note that this upper bound has also been obtained by using a computer
calculation with the computer algebra system {\sf GAP} \cite{gap}; see Table
\ref{table:orth}.

The $3$-transpositions in a parabolic subgroup of $G=\langle D\rangle$
stabilizing an isotropic point of $M$ generate a subgroup of $G$ which, up to
its center, is isomorphic to $3^5\split\Sym_6\cong 3^5\split W(A_5)$. The
$3$-transposition subalgebra generated by the $3$-transpositions in this
subgroup is, modulo its vanishing ideal, isomorphic to the Lie algebra of type
$^2A_5(2)$ modulo its center; see \PropRef{prop:weyl}. This latter algebra has
dimension $34$.

The above implies that $\OA$ has dimension $34$ and it is isomorphic to the Lie
algebra of type $^2A_5(2)$ modulo its center.
\end{proof}

\begin{lemma}\label{O6plus}
Suppose $D$ is one of the two classes of reflections in $\O_6^+(3)$. Then the
related algebra $\OA$ is isomorphic to a Lie algebra of type $D_4(2)$ modulo its
center.
\end{lemma}

\begin{proof}
The existence of a $3$-transposition subgroup of $G=\langle D\rangle$ isomorphic
to the group $3^4\split W(D_4)$ (as found inside a parabolic subgroup of $G$)
implies that the dimension of $\OA$ is at least $26$. Indeed, Lemma
\ref{lem:subalgebra} and \PropRef{prop:weyl} imply that $\OA$ contains a
subalgebra having a quotient isomorphic to the $26$-dimensional simple Lie
algebra of type $D_4$.

On the other hand, $G$ is a subgroup of $F_4(2)$ generated by long root
elements; see \cite{atlas}. These root elements are the $3$-transpositions in
$G$. The Lie algebra of type $F_4(2)$ is not simple. It has a $G$-invariant
ideal of dimension $26$ generated by the short roots. The quotient
algebra is a simple Lie algebra of type $D_4$. Hence, the algebra $\OA$ is
isomorphic to this latter algebra. 
\end{proof}

\begin{lemma}
Suppose $D$ is the class of $+$-reflections  in  $^+\Omega^+_7(3)$ or the unique
class of $3$-transpositions in $\Fi{22}$. Then $\OA$ is the Lie algebra of type
$^2E_6(2)$.
\end{lemma}

\begin{proof}
We note that $^+\Omega^+_7(3)$ is a $3$-transposition subgroup of $\Fi{22}$; see
\cite{asch,atlas}. Inside the group $^+\Omega^+_7(3)$, we find that the
reflections in $D$ that stabilize a fixed isotropic vector in the natural
orthogonal module for $^+\Omega^+_7(3)$ form a Fischer space isomorphic to that
of $3^6\split W(E_6)$. This shows that (in both cases) the algebra $\OA$
contains a subalgebra which has a quotient isomorphic to the Lie algebra of type
$^2E_6(2)$ of dimension $78$. Thus, in both cases the dimension of $\OA$ is at
least $78$. If the dimension of $\OA$ happens to be equal to $78$, then we can
conclude that, in both cases, $\OA$ is isomorphic to the simple Lie algebra of
type $^2E_6(2)$.

This implies that we only have to show that the dimension of $\OA$ is at most
$78$ in the case of $\langle D\rangle=\Fi{22}$. (Indeed, by Lemma
\ref{lem:subalgebra}, the dimension of $\OA$ in the case of $\Fi{22}$ cannot be
smaller than the dimension in the case of $^+\Omega^+_7(3)$.)

This, however, follows from the results of \cite{cush} on the geometry of this
Fischer group. Indeed, the results presented there imply that the dimension of
$\OA$ is at most $78$.

Note that we also verified this result with the use of {\sf GAP} \cite{gap}; see
Table \ref{tables}.
\end{proof}

We remark that, in order to obtain the upper bounds for the dimension of $\OA$
in the above results, we could have used the fact that $\POmega^-_6(3)$ embeds
into $\PSU_6(2)$ and that both $^+\Omega^+_7(3)$ and $\Fi{22}$ embed in
$^2E_6(2)$. However, with our present approach, these embeddings become
consequences of the above results.

Indeed, we can use them to prove Corollary \ref{cor:embed}:

\begin{cor}
\begin{enumerate}
\item The group $\PSU_6(2)$ contains a subgroup isomorphic to $^+\Omega^-_6(3)$
generated by root elements.
\item The group $^2E_6(2)$ contains subgroups generated by root elements, which
are isomorphic to $\Fi{22}$ respectively $^+\Omega^+_7(3)$.
\end{enumerate}
\end{cor}

\begin{proof}
By the above results, we find that the group  $H=\PO^-_6(3)$ embeds into the
automorphism group of the unitary Lie algebra $\mathfrak{psu}_6(2)$. Moreover, by
construction, the $3$-transpositions of $G$ correspond to root elements in
$\mathfrak{psu}_6(2)$. However, this implies that $H$ embeds into $G=\SU_6(2)$.
Under this embedding, the $3$-transpositions of $H$ are root elements in $G$.
Indeed, as we have seen in Section \ref{sec:cartan-normalizer}, the
$3$-transpositions in the subgroup $3^5\split\Sym_6$ of $H$ correspond to root
elements in $\mathfrak{psu}_6(2)$; and this clearly implies (a).

Similarly, we find $\Fi{22}$ and $^+\Omega^+_7(3)$ to be subgroups of the
automorphism group of the $^2E_6(2)$ Lie algebra defined by their Fischer
spaces. Again, the $3$-transpositions correspond to root elements in the Lie
algebra. This proves (b).
\end{proof}

\begin{lemma}
Let $D$ be the unique class of $3$-transpositions in
${\Omega}_8^+(2)\split\Sym_3$. Then $\OA$ is isomorphic to the $26$-dimensional
Lie algebra of type $D_4(2)$.
\end{lemma}
  
\begin{proof}
The Fischer space on $D$ can be partitioned into three parts, $D_1$, $D_2$ and
$D_3$, each forming a Fischer space of type ${\Omega}_8^+(2)$, such that any
point in one part is collinear with all points in the other parts. This implies
that any vanishing set of $D_i$ of even size is also a vanishing set of $D$.

As the lines of the Fischer subspace $D_i$ are vanishing in $D_i$, we find that
they generate a subspace of $\F_2D_i$ of codimension $8$. Indeed, modulo this
subspace we obtain the natural embedding of the Fischer space in the orthogonal
space $\O_8^+(2)$; see \cite{hall}. However, this implies that the vanishing
sets of even size in $D_i$ generate a subspace of codimension $9$. As every
affine plane is a vanishing set of $D$ of odd size, we find that $\V(D)$ has
codimension at most $9+9+9-1=26$ in $\A(D)$.

Now, $G=\POmega_8^+(2)\split\Sym_3$ embeds in $F_4(2)$ in such a way that the
$3$-transpositions are long root elements. The Lie algebra of type $F_4(2)$ is
not simple. It has a $G$-invariant ideal of dimension $26$ generated by the
short root elements. The quotient algebra is a simple Lie algebra of type $D_4$.
See also \ref{O6plus}. This proves that the latter algebra is isomorphic to
$\OA$.
\end{proof}

We completed our proof of \ThmRef{thm:main}.

\section{Some computational results}
\label{sec:comp-results}

In this final section we present some computational results on algebras defined by $3$-transpositions. A few of these computations repeat the
results obtained above. We therefore want to stress that the proofs in the
preceding sections are entirely independent of these computations and do not
rely on computer.

We retain the notation of the previous sections. Thus, assume that $G$ is a
group generated by its class of $3$-transpositions $D$, and denote by $\A$ the
$3$-transposition algebra of $D$. The vanishing ideal is denoted by $\V$. In
addition, let $\IAff$ denote the ideal of $\A$ generated by the set $\Aff$ of
affine planes in the Fischer space $\Pi(D)$.

Using {\sf GAP} \cite{gap}, we computed the dimensions of two quotients of $\A$,
one being $\OA=\A/\V$ and the other $A/\IAff$. Note that the latter is the
largest quotient of $\A$ that is a Lie algebra.

To determine the dimension of $\OA$ we have used the following.

\begin{prop}
Suppose $D$ is finite. Then the dimension of $\OA$ equals the $\F_2$-rank of the
adjacency matrix of the collinearity graph of the Fischer space on $D$.
\end{prop}

\begin{proof}
Note that this matrix, viewed over $\F_2$, coincides with the Gram matrix of the
form $\form{\cdot}{\cdot}$. Hence the rank of the matrix is equal to the
codimension of the radical.
\end{proof}

\subsection{Unitary groups}

\begin{center}
\begin{tabular}{|c|c|c|c|}
\hline
Group $G$ & $|D|$ & $\dim\A/\IAff$ & $\dim\OA$\\
\hline
$\SU_2(2)$ & 3 & 3 & 2 \\
$\SU_3(2)$ & 9 & 8 & 8 \\
$\SU_4(2)$ & 45 & 30 & 14 \\
$\SU_5(2)$ & 165 & 45 & 24 \\
$\SU_6(2)$ & 693 & 78 & 34 \\
$\SU_7(2)$ & 2709 & 119 & 48 \\
$\SU_8(2)$ & 10965 & 176 & 62 \\
$\SU_9(2)$ & 43605 & 249 & 80 \\
$\SU_{10}(2)$ & 174933 & 340 & 98 \\
\hline
\end{tabular}
\end{center}

\subsection{Orthogonal groups over $\F_3$}
\label{table:orth}

\begin{center}
\begin{tabular}{|c|c|c|c|}
\hline
Group $G$ & $|D|$ & $\dim\A/\IAff$ & $\dim\OA$\\
\hline
$^+\Omega^+_3(3)$ & 3    & 3 & 0\\
$^+\Omega^-_3(3)$ & 6    & 3 & 2\\
$^+\Omega^+_4(3)$ & 12   & 12 & 2\\
$^+\Omega^-_4(3)$ & 15   & 15 & 4\\
$^+\Omega^+_5(3)$ & 36   & 36 & 6\\
$^+\Omega^-_5(3)$ & 45   & 30 & 14\\
$^+\Omega^+_6(3)$ & 117   & 52 & 26\\
$^+\Omega^-_6(3)$ & 126  & 56 & 34\\
$^+\Omega^+_7(3)$ & 351  & 78 & 78\\
$^+\Omega^-_7(3)$ & 378  & 0 & 104\\
$^+\Omega^+_8(3)$ & 1080 & 0 & 260\\
$^+\Omega^-_8(3)$ & 1107 & 0 & 286\\
$^+\Omega^+_9(3)$ & 3240  & 0 & 780\\
$^+\Omega^-_9(3)$ & 3321  & 0 & 860\\
$^+\Omega^+_{10}(3)$ & 9801 & 0 & 2420\\
$^+\Omega^-_{10}(3)$ & 9882 & 0 & 2500 \\
\hline
\end{tabular}
\end{center}

\subsection{Sporadic groups}
\label{tables}

\begin{center}
\begin{tabular}{|c|c|c|c|}
\hline
Group $G$ & $|D|$ & $\dim\A/\IAff$ & $\dim\OA$\\
\hline
$\POmega_8^+(2)\split \Sym_3$ & 360 & 52& 26\\
$\POmega_8^+(3)\split \Sym_3$ & 3240 & 0 & 782\\
$\Fi{22}$        &3510 &  78 & 78 \\
$\Fi{23}$        &31671 &  0 & 782 \\
$\Fi{24}$        &306936 &  0 & 3774 \\
\hline
\end{tabular}
\end{center}

\subsection{Some $3$-transposition groups with normal $3$-group}

We consider the examples {\bf PR9-16} from \cite{cuha}. To perform our
calculations, we made use of the presentations of these groups as given in
\cite{haso}.

We include $\dim\V$ in the table, to highlight that it is the same for certain
related examples. E.g., the two cases of {\bf PR11} both have $\dim\V=141$. Both
are extensions of $\SU_5(2)$, and there we also have $\dim\V=141$. Also, in the
case {\bf PR12} we see that for $I=0,1$ or $3$ we find that $\dim\V=28$.
Likewise, {\bf PR13} is an extension of $^+\Omega^-_5(3)$, and in each case we
find $\dim\V=31$.

Another interesting tidbit is that $\dim\OA$ is the same for $^+\Omega^+_9(3)$
and {\bf PR15}, and for $^+\Omega^-_{10}(3)$ and {\bf PR16}.

\begin{center}
\begin{tabular}{|l|l|c|c|c|c|}
\hline
Reference \cite{cuha}& Group $G$ & $|D|$ & $\dim\A/\IAff$ & $\dim\OA$ & $\dim\V$\\
\hline
{\bf PR9} with $\abs{I}=1$& $3^7\split W(E_7)$ & 189 & 133 & 132 & 57 \\
{\bf PR10} with $\abs{I}=1$& $3^8\split W(E_8)$ & 360 & 248 & 248 & 112 \\
{\bf PR11} with $\abs{I}=1$& $3^{10}\split  (2\times {\SU}_5(2))$ & 1485 & 0 & 1344 & 141 \\
{\bf PR11}  with $\abs{I}=2$& $3^{10+10}\split (2\times{\SU}_5(2))$ & 13365 & 0 & 13224 & 141 \\
{\bf PR12}  with $\abs{I}=1$& $3^8\split (2^{1+6}\split{\SU}_3(2)')$ & 324 & 0 & 296 & 28 \\
{\bf PR12}  with $\abs{I}=2$& $3^{8+8}\split (2^{1+6}\split {\SU}_3(2)')$ & 2916 & 0 & 2888 & 28 \\
{\bf PR13}  with $\abs{I}=1$& $(3^5\cdot 3^5)\split {}^+\Omega_5^-(3)$ & 405 & 0 & 374 & 31 \\
{\bf PR14}  with $\abs{I}=1$&$(3^6\cdot 3^6)\split (3\cdot {}^+\Omega_6^-(3))$ & 1134 & 0 & 1042 & 92 \\

{\bf PR15}  with $\abs{I}=1$& $3^7\cdot {}^+\Omega_7^-(3)$ & 1134 & 0 & 860 & 274\\
{\bf PR16}  with $\abs{I}=1$& $3^8\cdot {}^+\Omega_8^-(3)$ & 3321 & 0 & 2500 & 821\\
\hline
\end{tabular}
\end{center}

\begin{remark}
The data in the above tables leads to some further interesting questions. For
example, the $3$-transposition algebra $\OA$ of $\POmega_8^+(3)\split\Sym_3$ has
the same dimension as the one associated to $\Fi{23}$. This implies that they
are the same. Can this observation lead to a new construction of the group
$\Fi{23}$ as a (sub)group of the automorphism group of the $3$-transposition
algebra associated to $\POmega_8^+(3)\split \Sym_3$? Is $\Fi{23}$ the full group
of automorphisms of this algebra?
\end{remark}

\end{document}